\documentclass{amsart}
\usepackage{amsmath,amsthm,amssymb,amsfonts}
\usepackage[colorlinks=true]{hyperref}
\usepackage{graphicx}
\usepackage{tikz}
\hypersetup{urlcolor=blue, citecolor=blue}

\def\doi#1{   {\href{http://dx.doi.org/#1}
   {{\mdseries\ttfamily DOI}}}}

\usepackage{graphics}

\newcommand{\gm}{\mathfrak{g}} 

\newcommand{\al}{\alpha}    \newcommand{\be}{\beta}
\newcommand{\de}{\delta}    
  \newcommand{\ep}{\varepsilon}
  \newcommand{\La}{\Lambda}
    \newcommand{\la}{\lambda}
    
    \newcommand{\Om}{\Omega}
\newcommand{\ga}{\gamma}    
\newcommand{\R}{\mathbb{R}}
\newcommand{\N}{\mathbb{N}}

\newcommand{\pt}{\partial_t}\newcommand{\pa}{\partial}
\newcommand{\les}{{\lesssim}}
\newcommand{\beeq}{\begin{equation}}\newcommand{\eneq}{\end{equation}}

\newcommand{\Sp}{{\mathbb S}}\def\CO{\mathcal {O}}

\newcommand{\supp}{\text{supp}}

\newenvironment{prf}{\noindent {\bf Proof.} }{\endprf\par}
\def \endprf{\hfill  {\vrule height6pt width6pt depth0pt}\medskip}

\def\<{\langle}             \def\>{\rangle}
\def\({\left(}                 \def\){\right)}
\numberwithin{equation}{section}
\newtheorem{thm}{Theorem}[section]

 \newtheorem{defn}[thm]{Definition}
 \newtheorem{rem}[thm]{Remark}

\title[Wave equations on Kerr black hole backgrounds]
{Blow-up for semilinear wave equations on Kerr black hole backgrounds}

\author{Mengyun Liu}
\address{Department of Mathematics\\Zhejiang Sci-Tech University\\Hangzhou 310018, P. R. China }
\email{mengyunliu@zstu.edu.cn}

\author{Chengbo Wang}
\address{School of Mathematical Sciences\\ Zhejiang University\\Hangzhou 310058, P. R. China}\email{wangcbo@zju.edu.cn }



\date{\today}

\thanks{The first author was supported by NSFC 12101558 and NSF of Zhejiang province LQ22A010016.
The second author was supported by
 NSFC 11971428  and  NSFC 12141102. }

\begin{document}
\maketitle

\begin{abstract}
We examine solutions to semilinear wave equations on black hole backgrounds and give a proof of an analog of the  blow up part of the John theorem,
with $F_p(u)=|u|^{p}$,
 on the Schwarzschild and Kerr black hole backgrounds. Concerning the case of Schwarzschild,
we construct a class of small data, so that the solution blows up along the outgoing null cone, which applies
for both $F_p(u)=|u|^{p}$ and the focusing nonlinearity $F_p(u)=|u|^{p-1}u$.
The proof suggests that the black hole does not have any essential influence on the formation of singularity, in the region away from the Cauchy horizon $r=r_-$ or
the singularity $r=0$.
Our approach is also robust enough to be adapted for
general asymptotically flat space-time manifolds,
possibly exterior to a compact domain, with spatial dimension $n\ge 2$.
Typical examples include exterior domains, asymptotically Euclidean spaces,
Reissner-N\"ordstr\"om space-times, 
and Kerr-Newman  space-times.
\end{abstract}


\section{Introduction}
The purpose of this article is to 
 investigate the blow up phenomenon for solutions to nonlinear wave equations with small initial data on 
the Schwarzschild and Kerr black hole backgrounds.
 In particular, we establish the blow up result for
solutions to a class of semilinear wave equations with power-type
nonlinearities with power less than a certain critical power.  This
critical power, $1+\sqrt{2}$, is the same as that on
$(1+3)$-dimensional Minkowski space-time, which was known from the seminal work of John \cite{John79}.

Let us first recall the Kerr metric,
which is a stationary
solution to Einstein's vacuum equations for which we have parameters $a$ and $M$.  In Boyer-Lindquist coordinates, with  $0\le a< M$, it is given by
\beeq\label{metric-Kerr}
ds^2 = g_{tt}\,dt^2 +2 g_{t\phi} dt\,d\phi + g_{rr}\,dr^2 
 + g_{\theta\theta}\,d\theta^2+g_{\phi\phi}\,d\phi^2\ ,
\eneq
where 
$t \in \R$, $r > 0$, $(\theta, \phi)$ are the spherical
coordinates on $\Sp^2$ and
\[
g_{tt}=-\left(1-\frac{2Mr}{\rho^2}\right), \qquad
g_{t\phi}=-\frac{2Ma r\sin^2\theta}{\rho^2}, \qquad
g_{rr}=\frac{\rho^2}{\Delta},
\]
\[g_{\theta\theta}={\rho^2},\qquad 
 g_{\phi\phi}=\frac{(r^2+a^2)^2-a^2\Delta
  \sin^2\theta}{\rho^2}\sin^2\theta, 
\]
with
\[
\Delta=r^2-2Mr+a^2, \qquad \rho^2=r^2+a^2\cos^2\theta.
\]
Here $M$ represents the mass of the black hole and $aM$ its angular
momentum.
The Schwarzschild space-time is the static solution
corresponding to $a=0<M$.  And the Minkowski space-time is the trivial
solution for which we have $a=0$ and $M=0$.

For $0<a<M$, there are  apparent
 singularities at the
 two roots of 
$\Delta =0$,
$$r_{\pm}=M\pm\sqrt{M^2-a^2}\ .$$
It turns out
they are merely coordinate singularities.
The  region $r\in (r_-, r_+)$ is a black hole region for the 
Kerr space-time. The surface,
 $r=r_+$, is the outer event horizon, which is the past null boundary of the black hole. The surface,
 $r=r_-$, is the Cauchy horizon, which is the future null boundary of the black hole and it is the boundary of the maximal region so that we could have well-posed theory for the Cauchy problem in the Kerr black hole. 
 For a further discussion of the nature of $r_\pm$, which is
not relevant for our results, we refer the reader to, e.g.,
 \cite{MR0424186}, \cite{MR700826}.

Let $p>1$, we will consider the evolution of the nonlinear
waves on Kerr black hole backgrounds,
\begin{equation}\label{nlw}
\Box_K u + F_p(u)=0, \qquad u |_{\tilde v=0} = u_0, \qquad \tilde T u |_{\tilde v=0} = u_1\ .
\end{equation}
Here $\Box_K$ denotes the d'Alembertian in the Kerr metric, 
 the coordinate $\tilde v$
is chosen so that the slice $\tilde v=0$ is uniformly space-like, and
$\tilde T$ is
a smooth, uniformly time-like vector field.
For clarity, we require $\tilde v=t$ 
and $\tilde T=\pa_t$ near spatial infinity, $r\gg M$.
See, e.g,
\cite{LMSTW} for more details of the coordinates.

We shall assume that the nonlinear term behaves like $|u|^p$ \beeq\label{key}
 \sum_{0\leq j\leq [p]}  |u|^j |\partial_u^j F_p(u)| \leq C |u|^{p} \ .
\eneq
Typical examples include $F_p(u)=\pm |u|^p$ and $\pm |u|^{p-1}u$.

\subsection{Results for Minkowski space-time}
In Minkowski space-time  $\R^{n, 1}$ with $n\ge 2$ and  $F_p(u)=|u|^p$, the problem with small data is known as the John problem or the Strauss conjecture.

In the first work in this area,
 F.  John   \cite{John79} proved that the power $1+\sqrt{2}$ is critical in the sense that, we have global existence for small compactly supported smooth data when $p>1+\sqrt{2}$, while blow-up could occur for arbitrarily small data when  $p<1+\sqrt{2}$. 
 
 Shortly afterward, 
 Strauss
  \cite{Strauss81} gave the conjecture that the critical power $p_c(n)$ is the positive
root of the quadratic equation
$$(n-1)p^2-(n+1)p-2=0\ ,$$ with $p_c(3)=1+\sqrt{2}$ and $p_c(4)=2$.
  The existence portion of the conjecture 
  (for supercritical, but sub-conformal powers,
   $p\in (p_c, 1+4/(n-1)]$)
  was verified in 
  \cite{Glassey81ex}
($n=2$), \cite{Zhou95} ($n=4$), \cite{LdSo96} ($n\le 8$ or radial data), and
\cite{GLS97, Ta01-2} (all $n\ge 3$).
Concerning the current status of the art for the super-conformal powers, we refer the readers to the recent work
\cite{SW22}.

On the other hand, when $p\in (1, p_c(n))$,
the blow up results,
as well as the following upper bound estimates of the lifespan, denoted by $T_*$,
\beeq\label{eq-upperbd}
T_*\le C \ep^{\frac{2p(p-1)}{(n - 1)p^2 - (n + 1)p - 2 }}\ ,
\eneq
 have been obtained
through the works \cite{Glassey81bu}, \cite{Sideris84}, for certain compactly supported smooth initial data of the form
$(u_0, u_1)=(\ep\phi, \ep\psi)$, with $0<\ep\ll 1$.
The upper bound is also sharp in general, at least for $n\ge 3$ or $p>2/(n-1)$, for which we refer
 \cite{LLWW21} and references therein for more details.
See \cite{Scha85}, 
 \cite{LdSo96},
\cite{YorZh06} and \cite{Zh07} for the corresponding results for the critical case.


\subsection{Global theory}
Turning back to the Kerr space-time.
For sufficiently small data with sufficient fast decay, the global well-posedness
 has been established when $p> 1+\sqrt{2}$, for the Kerr space-time with small angular momentum $a\ll M$.


In the sample case of the
 Schwarzschild  black hole backgrounds, 
 the global existence of the Cauchy problem \eqref{nlw}, in the region $r>r_0$ (with fixed $r_0\in (r_-, r_+]=(0, 2M]$), has been obtained in \cite{LMSTW}, for any $p>1+\sqrt{2}$, which is an analog of the existence part of the John theorem
  \cite{John79} from the Minkowski space-time to the Schwarzschild space-time. Previous works include
   \cite{DaRod05} (for
the general  Reissner-Nordstr\"om space-time, with
  $p>4$, $r>2M$ and radial data), \cite{BlueSter06} ($p>3$,  $r>2M$).

Concerning the general Kerr space-time, with sufficiently small angular momentum $a\ll M$, the global result (with $p>1+\sqrt{2}$) is obtained in \cite{LMSTW}, under the assumption that the initial data are compactly supported. The technical support assumption was removed later in
   \cite{MW17}.


\subsection{Long time existence}
For the case with subcritical or critical powers, the (long time) well-posed theory has also been developed, in \cite{W17},
for the Schwarzschild space-time $a=0$.
More precisely, it is shown that, there exist
$R\gg 3M$ and
$c=c(p)>0$ so that, we have 
well-posedness, in $C_{\tilde v}H^3\cap
C_{\tilde v}^1H^2$, in the domain $(\tilde v, r, \omega)\in [0, T_*)\times (r_0,\infty)\times\Sp^2$,
with
\beeq\label{eq-lowerbd}
T_*\ge T_\ep=
\begin{cases}
c \ep^{-\frac{p(p-1)}{1+2p-p^2}}\ ,\ 2\leq p< 1+\sqrt{2}\ ,\\
e^{c \ep^{-2(p-1)}}\ , ~~~~~~~~~~~~~~~\ \ p=1+\sqrt{2}\ ,\ 
\end{cases}
\eneq
whenever the initial data $(u_0, u_1)\in H^3\times H^2$ satisfies
\beeq
\label{initial}
\|(u_0, u_1)\|_{X}:=
\|\phi_R (u_0, u_1)\|_{H^3\times H^2}+
\sum_{|\ga|\le 2} \|\psi_R (\nabla, \Omega)^{\ga}(\nabla u_0, u_1)\|_{L^2\cap\dot{H}^{s_d-1}}\le \ \ep\ll 1\ .
\eneq
Here
$$s_d=\frac{1}{2}-\frac{1}{p}\ ,$$
$\Omega$ stands for the rotational vector fields, $\psi_R(r)=\psi(r/R)$, $\phi_R(r)=1-\psi_{2R}$, 
with  smooth $\psi(r)=0$ for $r<1$ and $\psi=0$ for $r>2$.
We remark here that we have stated the result with a slightly weaker assumption on the initial data, \eqref{initial}, than what stated in \cite[Theorem 6.1]{W17}. As is clear from the proof, this assumption is sufficient for the same result. As we shall see, the reason for us to state in the current form is for the convenience of the discussion of the blow up theory.

As a matter of fact, it is interesting to recall that, the exponent occurred in the lifespan is related with $s_d$,
\beeq
\frac{2p(p-1)}{(n-1)p^2-(n+1)p-2}=\frac{1}{s_c-s_d} \ ,\eneq
where $s_c=n/2-2/(p-1)$ is critical Sobolev regularity with respect to the scaling.

For the case $0<a\ll M$, a weaker result for $2\le p\le 1+\sqrt{2}$
 was also available
in \cite[Theorem 7.1]{W17}
for compactly supported data (and shorter lower bound of the lifespan for the critical case).
These technical restrictions was relaxed later in
\cite[Corollary 4.2]{MR4112037} for $2<p\le 1+\sqrt{2}$.
In particular, when $2<p\le 1+\sqrt{2}$,
we have \eqref{eq-lowerbd} for any data satisfying \eqref{initial}.


\subsection{Upper bound of the lifespan}

At first, for clarity, we define the lifespan of the problem
\eqref{nlw}, with size $\ep$, as follows.
\begin{defn}[Lifespan]
Let  $T_*(u_0, u_1)$ be the supremum of $T>0$ so that 
the problem \eqref{nlw} is well-posed in the domain
 $(\tilde v, r, \omega)\in [0, T)\times (r_0,\infty)\times\Sp^2$,
for given data $(u_0, u_1)$ satisfying \eqref{initial}.
The  infimum of all such time of existence, \beeq
T_*(\ep)=\inf_{\|(u_0, u_1)\|_X\le \ep}T_*(u_0, u_1)\ ,
\eneq
is called to be the lifespan of the problem
with size $\ep>0$.
\end{defn}

In view of the well-posed theory, we have
\beeq\label{eq-lowerbd0}
T_*(\ep)\ge 
\left\{\begin{array}{ll}
c\ep^{-2}\ ,& p=2, a=0\ ,\\
c\ep^{-\frac{p(p-1)}{1+2p-p^2}}\ ,& 2< p< 1+\sqrt{2}, 0\le a\ll M\ ,\\
e^{c\ep^{-2(p-1)}}\ , & \ p=1+\sqrt{2},
 0\le a\ll M
\ . 
\end{array}\right.
\eneq
Then it is natural to ask what is the sharp estimate for the lifespan.

Observing that the lower bound, \eqref{eq-lowerbd0}, for the subcritical case, $2\le p<p_c$, agrees with the upper bound 
on the Minkowski space-time, \eqref{eq-upperbd}, which suggest that it should be sharp.

Unlike the well-posed theory,
blow up results are available only for the  Schwarzschild black hole space-time, $a=0< M$, with $F_p(u)=|u|^p$, in the black hole exterior domain $r>r_+=2M$, where
\beeq\label{metric-Sch}
ds^2=-(1-\frac{2M}{r})dt^2+(1-\frac{2M}{r})^{-1}dr^2+r^2d\omega^2\ .\eneq
See \cite{CaGe06} for a weaker blow up result for $1<p< 1+\sqrt{2}$ with data of the form
$(\ep\phi(r-\ep^{-N}), \ep\psi(r-\ep^{-N}))$ for certain $N\gg 1$, and
  \cite{LinLaiMing19} for blow up result for $1<p\le 2$ with data of the form $(\ep\phi, \ep\psi)$. 
Recently, Lai and Zhou \cite{LZ22} obtained the blow up result, as well as the expected upper bound of the lifespan, in the $(t,r)$ coordinates,
\eqref{eq-upperbd},  for $2\le p< 1+\sqrt{2}$ with data of the form $(\ep\phi(r), \ep\psi(r))$, which are compactly supported, away from the event horizon.

\subsection{Main results}
Before discussing our results. It is interesting to compare what is known in the Schwarzschild space-time.
The well-posed theory developed in \cite{W17} tells us that the solutions remain bounded in the region
$$\{(\tilde v, r, \omega): \tilde v\in [0, T_\ep], r>r_0, \omega\in\Sp^2\}\ ,$$
where $T_\ep$ is that appeared in
\eqref{eq-lowerbd}. Moreover, as long as the solution has bounded $X$ norm on the Cauchy surface of constant $\tilde v$, the solution could be extended for larger $\tilde v$.
In addition, the initial surface $\tilde v=0$ could be replaced by the Cauchy surface of the Schwarzschild space-time, corresponding to $t=0$. The well-posedness from $t=0$ to $\tilde v=0$ follows from a standard energy argument.

On the other hand, the blow up result of
\cite{LZ22} tells us that, in the $(t,r)$ coordinates,
we have an upper bound of the lifespan,
\eqref{eq-upperbd}.

Noticing, however, that in the part with $r<R$, where $R\gg 1$ is appeared in the $X$ norm, 
the region $$\{(t,r,\omega): 
t\in [0, \infty), 2M<r<R, \tilde v\le T_\ep, \  \omega\in \Sp^2\}
$$ lies in the domain of existence.
Based on this observation, we
see that, for $t=T_*$, the solution shall blow up in the part $r>R$, which is close to the spatial infinity.

In light of this comparison,
we would like to construct a class of data, small in $X$ norm, so that the solution blows up along the outgoing null direction.

\subsubsection{Schwarzschild space-time}
At first, concerning on the sample case of Schwarzschild space-time ($a=0< M$),
by finding a sufficiently small $\ep_0\in (0, 1)$,
%
we construct a class of outgoing radial data $(u_0, u_1)$ of the form
\beeq\label{eq-data-a=0} u_0(r)=
\ep_0^{-\theta_0}
r^{-\al-1},
u_1(r)=
\ep_0^{-\theta_1}
r^{-\al-2}, 
r\in [\ep^{-N}, 10\ep^{-N}]\ , 
\eneq
where $0<\theta_0<\theta_1<1$,
$$
\al=\frac{2}{p-1}-1-\frac \mu {p(p-1)},\ 
N=N_\mu=\frac{p(p-1)}{1+2p-p^2-\mu}>
N_0=\frac{p(p-1)}{1+2p-p^2}\ ,
$$
for any $\ep, \mu\in (0,\ep_0)$.
It is clear that such data \eqref{eq-data-a=0} could be extended to be $C_0^\infty$ functions so that
\beeq\|(u_0, u_1)\|_{X}\le \ep_0^{-1}\ep, \forall \mu, \ep<\ep_0\ .\eneq
With such data, we are able to prove that in the domain of determination of
$r\in(\ep^{-N}, 10\ep^{-N})$, the solution remains positive, 
concentrate along the outgoing
 light cone
$$r^*-t=2\ep^{-N}\ ,
$$
and blows up in finite time:
\beeq  T_*(u_0, u_1)\le 3\ep^{-\frac{p(p-1)}{1+2p-p^2-\mu}}, \forall \ep, \mu\in (0,\ep_0)\ .\eneq
Here,
$r^*=r+2M\ln (r-2M)$ 
is the tortoise coordinate so that
$dr^*=(1-2M/r)^{-1}dr$.
In conclusion, letting $\mu$ goes to zero, we arrive at our first blow up result, for the Schwarzschild space-time.

\begin{thm}
\label{th3}
Let $1<p<1+\sqrt{2}$, $F_p(u)=u^p$ when $u>0$. Consider the Cauchy problem \eqref{nlw} posed on Schwarzschild space-time ($a=0< M$). Then
there exist $C, \ep_0>0$ so that, 
\beeq\label{eq-upperbd0}
T_*\le C\ep^{-\frac{p(p-1)}{1+2p-p^2}}\ ,
\eneq
for any $\ep\in (0,\ep_0)$.
Here, the constant could be set to be
$$C=3\ep_0^{-\frac{p(p-1)}{1+2p-p^2}}\ .$$
\end{thm}
Notice that one of the novelties of our approach is that our result applies also for the focusing nonlinearity $F_p(u)=|u|^{p-1}u$.
Also, it illustrates the blow up phenomenon vividly, which blows up
along the outgoing
 light cone.

\subsubsection{Kerr space-time}
Inspired by the result of the Schwarzschild space-time, we turn to the more general Kerr space-time.
We will choose data of the similar form, although the previous approach does not apply in the case of $a>0$, which rely heavily on the spherically symmetric nature of the  Schwarzschild space-time. Instead, we use the functional method, which has been widely used in proving blow up results, for the Strauss conjecture on  backgrounds other than the Kerr space-time.

Remarkably, thanks to the structural condition
\beeq\label{eq-stru-Kerr}
\pa_\be(g^{1/2}g^{0\be})=0\ ,\eneq
where
$g^{\al\be}$ is the inverse of $g_{\al\be}$ and $g=-\det (g_{\al\be})$,
it turns out that the approach for the Kerr space-time is much simpler than that for the
Schwarzschild space-time, with the price that it works only for the nonlinearity $F_p(u)=|u|^p$, as well as the lack of the precise knowledge on the blow up mechanism.
\begin{thm}\label{thm-Kerr}
Let $1<p<1+\sqrt{2}$, $F_p(u)=|u|^p$. Consider the Cauchy problem \eqref{nlw} posed on Kerr space-time ($0\le a\le M$). Then
there exist $C, \ep_0>0$ so that, for any $\ep\in (0,\ep_0)$, we have
\eqref{eq-upperbd0},
for a class of small data in $C^{\infty}_0(r> 4R)$ satisfying
\eqref{initial}. In particular, the lower bound 
\eqref{eq-lowerbd0}
for $2< p< 1+\sqrt{2}$ is sharp in general.
\end{thm}

\subsection{Further generalization: asymptotically flat space-time manifolds}
Equipped with the blow up results for the Kerr space-time, it is natural to ask to what extent it could be generalized for 
asymptotically flat space-time.

Let $n\ge 2$,
we shall work on the $(1+n)$-dimensional asymptotically flat Lorentzian space-time $(\mathcal M,\gm)$.
The smooth Lorentzian metric, $\gm=  g_{\al\be} (t,x) \, d x^\al \, d x^\be$,
with inverse $g^{\al\be}$, $(t,x)=(x^0, x^1, \cdots, x^n)$,  is assumed to be asymptotically flat. That is, there exist $R>0$ and $\de>0$ such that
\begin{equation}\tag{H1} \label{eq-H1}
\partial_{tx}^\ga  (g_{\al\be}(t,x)-m_{\al\be})=\CO(|x|^{-|\ga|-\de}), |x|> R, \forall \ga\in\N^{1+n}\ ,
\end{equation}
where $(m_{\al\be})=Diag(-1,1,1,\cdots,1)$ is the standard Minkowski metric.

Moreover, the  space-time $(\mathcal M,\gm)$ is said to be spherically symmetric, if, in the polar coordinate
$x=r\omega$
with $\omega\in\Sp^{n-1}$, it could be written in the form
\beeq\label{metric-sph}\gm=\gm_0=g_{tt}(t,r)dt^2+2g_{tr}(t,r)dtdr+g_{rr}(t,r)dr^2+ r^{2}d\omega^2, \omega\in \Sp^{n-1}\ .\eneq
In this setting, the  asymptotically flat assumption \eqref{eq-H1} with parameter $\de=\de_0$ could be rephrased as
\begin{equation}\tag{H2} \label{eq-H2}
\partial_{t}^j\partial_{r}^k  (g_{\al\be}(t,r)-m_{\al\be})=\CO(r^{-\de_0-j-k}), r> R, \forall j,k\ .
\end{equation}

 Given $p>1$ and a Cauchy surface $\Sigma\subset \mathcal M$,  together with a uniformly timelike vector $\mathcal N$,
we shall consider the following Cauchy problem
\beeq
\label{eq-AF}
\begin{cases}
\Box_{\gm}u:=\nabla^\al\pa_\al u
=-F_p(u),
(t, x)\in \mathcal M\\
u|_\Sigma=u_0\ , \mathcal N u|_{\Sigma}=u_1\ .
\end{cases}
\eneq
For simplicity,
we assume $$\Sigma_R:=\Sigma\cap \{r>R\}=
\{(0,r,\omega): r>R, \omega\in \Sp^{n-1}\},\ 
\mathcal N=\pt, \forall r>R\ .
$$

Noticing that the setting of
 asymptotically flat space-time manifolds includes a lot of important space-time manifolds in general relativity, analysis and geometry.
Besides what we have mentioned for the Schwarzschild and Kerr space-time, the Strauss conjecture has also been investigated for
many nontrivial backgrounds.
These include 
\cite{DMSZ}, 
\cite{HMSSZ},
\cite{SSW12},
which examine global existence for similar equations exterior to nontrapping obstacles, with $2\le n\le 4$. See
 \cite{LLWW21} and references therein
 for blow up results.
In the case of nontrapping asymptotically Euclidean  manifolds, the global results were obtained in 
\cite{SW10}, \cite{WaYu11} for $n = 3, 4$, while the blow up result was known from
\cite{WaYo18-1pub},
\cite{LW19-p}, for certain exponentially perturbation of the 
spherically symmetric time-independent metric
\eqref{eq-H2}.
It also includes the 
Reissner-Nordstr\"om space-time,
Kerr-Newman  space-time, 
and time-dependent asymptotically flat perturbations of the Minkowski space-time, for which 
the  well-posed theory has been developed in the general context of  asymptotically flat space-time, in
   \cite{MW17}.
\cite{W17} and
\cite{MR4112037}.


A standard energy argument shows that there exists a local in time energy solutions to this problem, in the future domain of determination of $\Sigma_R$, denoted by $\mathcal M_R$. With $t$ as the parameter of foliation, we have $\mathcal M_R=\cup_{t\ge 0} \mathcal M^t_R$, and
$$u\in C([0,T]; H^1(\mathcal M^t_R))\cap C^1([0,T]; L^2(\mathcal M^t_R))\cap 
C^2([0,T]; H^{-1}(\mathcal M^t_R))\ .$$
The lifespan $T_*$, for the solution with data $(u_0, u_1)$, is the supremum of the time of existence $T$.

The  well-posed theory developed in \cite{W17} and
\cite{MR4112037}  applies to the current setting, provided that
$2\le p\le p_c(n)$, $n=3,4$ and
$$\gm=\gm_0+\gm_1-m$$
where
$\gm_0$ is spherically symmetric, \eqref{eq-H2}, with $\de_0>0$ and
$\gm_1$ 
satisfies \eqref{eq-H1} with $\de>1$.
More precisely, there exists $R_0\gg 1$,  so that,
for any data satisfying
\beeq
\label{eq-initial}
\sum_{|\ga|\le k} \| (\nabla, \Omega)^{\ga}(\nabla u_0, u_1)\|_{L^2\cap\dot{H}^{s_d-1}}\le \ \ep\ , 
\ s_d=\frac{1}{2}-\frac{1}{p}\ ,
\eneq
with $k=2$,
we have 
\beeq\label{eq-lowerbd2}
T_*\ge T_\ep=
\begin{cases}
c\ep^{-\frac{p(p-1)}{1+2p-p^2}}\ ,\ 2\leq p< 1+\sqrt{2}, n=3\ ,\\
e^{c\ep^{-2(p-1)}}\ , ~~~~~~~~~~~~~~~\ \ p=p_c(n), n=3,4\ ,\ 
\end{cases}
\eneq
in the region $\mathcal{M}_{R_0}$.

Our blow up result on asymptotically flat manifolds 
is given in the following:

\begin{thm}\label{thm-AF}
Let $n\ge 2$ and $1<p<p_c(n)$. Consider the Cauchy problem \eqref{eq-AF} posed on asymptotically flat manifolds \eqref{eq-H1}, 
with
\eqref{eq-stru-Kerr} and
 $(u_0, u_1)\in C^{\infty}_0(r> R)$.
For any $k>0$, there exists $\ep_0>0$ such that for any $\ep\in (0,\ep_0)$, there exists $(u^\ep_0, u^\ep_1)\in C^{\infty}_0(|x|> R)$ satisfying
\eqref{eq-initial},
such that 
the solution to the Cauchy problem \eqref{eq-AF} blows up in finite time. Moreover, the lifespan of the solution 
satisfies 
$$T_{*}\le 3(\ep_0\ep)^{-\frac{2p(p-1)}{2+(n+1)p-(n-1) p^2}}\ .$$
In addition,
when the
space-time manifold is spherically symmetric
\eqref{eq-H2} (without assuming the structural condition 
\eqref{eq-stru-Kerr}), the same result applies for
$F_p(u)=u^p$ when $u>0$.
\end{thm}

\subsection{Outline}
Our paper is organized as follows. In  the next section, for the sample case of the Schwarzschild space-time, we present our proof
of Theorem \ref{th3}. In this case, it is well known that, for radial solutions, the problem 
\eqref{nlw} could be transformed into a $1+1$ dimensional wave equation, with a short range potential,
see \eqref{eq-Sch-1d-2}. With the help of this formulation and iteration, the initial data could be exploited to obtain certain lower bound, which grows fast along the outgoing light cone.
 To illustrate blow up phenomenon, we employ two different approaches. One 
approach is to derive an ordinary differential inequality for certain functional, localized near the light cone,
see Section \ref{sec-2.4}.
Another approach is to obtain sharp estimate on the constant, which reveals blow up for the solution instead of a functional,
see Section \ref{sec-2.5}. We think either approach has their own advantage, and we choose to present both of them.

In Section \ref{sec-Ker}, we consider the similar problem for the Kerr space-time. In contrast to the Schwarzschild space-time, which is both static and spherically symmetric, the idea of the proof for the Schwarzschild space-time does not apply in the general case of $a>0$. Despite that, we are able to prove similar result, when $F_p(u)=|u|^p$, by an indirect approach. More precisely,
thanks to the structural condition \eqref{eq-stru-Kerr},
 we introduce certain functional and establish the blow up result by
 deriving an ordinary differential inequality.

In the remaining sections, we try to generalize the results to 
general asymptotically flat space-time.
It turns out that the stationary assumption is not necessary for similar result, 
see Section \ref{sec-AF} (under the structural condition \eqref{eq-stru-Kerr}).

In the last  Section \ref{sec-AF-radial}, for the 
spherically symmetric asymptotically flat space-time (without the condition \eqref{eq-stru-Kerr}),
we develop a similar blow up results for the Schwarzschild space-time. In this general setting, we employ the double null frame, in which, we could integrate the equation along the characteristics, to obtain similar lower bound and blow up phenomenon.

 \section{Schwarzschild space-time}\label{sec-Sch}

As we have mentioned in the introduction, we will choose
data of the form
\eqref{eq-data-a=0}. 
However, as we shall see, without loss of generality, it is more convenient to assume the data keeping such form near spatial infinity.
That is,
with $\ep_0, \theta_0, \theta_1\in (0, 1)$ to be specified, for any $\mu,\ep\in (0,\ep_0)$,  we set
\beeq\label{eq-data-a=0-1} u_0(r)=
\ep_0^{-\theta_0}
r^{-\al-1},
u_1(r)=
\ep_0^{-\theta_1}
r^{-\al-2}, 
 r\ge\ep^{-N}\ , 
\eneq
where
$$
\al=\frac{2}{p-1}-1-\frac \mu {p(p-1)},\ 
N=N_\mu=\frac{p(p-1)}{1+2p-p^2-\mu}>
N_0=\frac{p(p-1)}{1+2p-p^2}\ .
$$ Here we ask
$\ep_0<1+2p-p^2$ and
$\ep_0^{-N}>R\gg 3M$, so that $\al>1/p$. Noticing that $N(\al-1/p)=1$, which ensures
they
 could be extended to be $C^\infty$ functions so that
\beeq\|(u_0, u_1)\|_{X}\le \ep_0^{-1}\ep, \forall \mu, \ep<\ep_0\ .\eneq

\subsection{Transformation} 
For the Schwarzschild space-time,
in the 
domain of the future determination of
$\{(t,r): t=0, r\ge\ep^{-N} \}$,
 it is well known that,  we could transform the
Cauchy problem \eqref{nlw} into a $1+1$ dimensional wave equation,
when the solution is spherically symmetric, $u=u(t,r)$.

To be more specific, let
$r^*=r+2M\ln (r-2M)$ 
be the tortoise coordinate so that
$dr^*=(1-2M/r)^{-1}dr$. In $(t,r^*,\theta,\phi)$
coordinates,
the Schwarzschild metric
\eqref{metric-Sch} is in the form
$$
ds^2=-\left(1-\frac{2M}{r}\right)(dt^2-(dr^*)^2)+r^2d\omega^2\ .
$$
Then
$g^{1/2}=\sqrt{-\det (g_{\al
\be})}=(1-\frac{2M}{r})r^2\sin\theta$, and
\begin{eqnarray*}\Box_\gm u(t,r)&=&
-g^{-1/2}\pa_t \left(1-\frac{2M}{r}\right)^{-1}g^{1/2}\pa_t u+
g^{-1/2}\pa_{r^*} \left(1-\frac{2M}{r}\right)^{-1}g^{1/2}\pa_{r^*}u\\
&=&\left(1-\frac{2M}{r}\right)^{-1}[-\pa_t^2 u+
r^{-2}\pa_{r^*} r^2\pa_{r^*}u]\\
&=&\left(1-\frac{2M}{r}\right)^{-1}[-\pa_t^2 u+
\pa_{r^*}^2 +2
\frac{\pa_{r^*}r}{r}
\pa_{r^*}u]\\
&=&\left(1-\frac{2M}{r}\right)^{-1}r^{-1}[-\pa_t^2 +
\pa_{r^*}^2 -
\frac{\pa^2_{r^*}r}{r}
](ru)\ .
\end{eqnarray*}

Let $W=ru$, noticing that
$$\frac{\pa^2_{r^*}r}{r}=\frac{2M}{r^3}\left(1-\frac{2M}{r}\right)\ ,$$
 the equation \eqref{nlw} is equivalent to  a $1+1$ dimensional wave equation with potential:
\beeq\label{eq-Sch-1d}
\left(\pa_t^2 -
\pa_{r^*}^2 +\frac{2M}{r^3}\left(1-\frac{2M}{r}\right)\right)W=
-(r-2M)\Box_\gm u
=(r-2M)F_p\left(\frac W r\right)\ .
\eneq

When $\ep_0$  is sufficiently small, we have
$r\le r^*\le 2 r$ for any $r\ge\ep_0^{-N_0}$.
Let $L(\mu):=T_*(u_0, u_1)\ge c\ep^{-N_0}$ be the lifespan
of the solution, in the domain of determination of 
$\{(t,r): t=0, r\ge \ep^{-N}\}$.
Under these assumptions,
we know that the region \beeq \label{eq-Sch-domain}\Omega:=\{(t,r^*): t\in [0, L), 
r^*\ge t+2\ep^{-N}
\} \ ,\eneq is a smaller region  than the domain of determination of 
$\{(t,r): t=0, r\ge \ep^{-N}\}$,
We will 
exploit the equation \eqref{eq-Sch-1d} to show blow up in the region $\Om$.

\subsection{Bootstrap: initial lower bound}
To estimate the size of $W$, we write
 the equation \eqref{eq-Sch-1d} in the following form
\beeq\label{eq-Sch-1d-2}
(\pa_t^2 -
\pa_{r^*}^2)W=
\left(1-\frac{2M}{r}\right)\left[rF_p\left(\frac W r\right)
- \frac{2M}{r^3}W\right]:=G(r^*, W)
\ .
\eneq

We claim that for any $(t, r^*)\in \Om$, we have 
\beeq\label{eq-Sch-lower-1}
W(t,r^*)\geq (4M)^{\frac{1}{p-1}}r^{\frac{p-4}{p-1}}>0\ ,
\eneq
which ensures
$$
rF_p\left(\frac W r\right)
\ge 2 \frac{2M}{r^3}W,\ 
G(r^*, W)>\frac{1}{2}\left(1-\frac{2M}{r}\right) rF_p\left(\frac W r\right)>
\frac{W^p}{3(r^*)^{p-1}} >0
\ .$$
In fact, for any $a\gg 1$, let 
$$X_a=\{\tau\in [0, L)\cap [0, a];
W(t, r^*)> (4M)^{\frac{1}{p-1}}r^{\frac{p-4}{p-1}}, \forall r^*\in [t+2\ep^{-N}
, a-t]
\ , \forall t\in [0, \tau]\}\ .$$
Then, by continuity, it is easy to see that $X_a$ is open and nonempty since 
$$W(0, r^*)=ru_0(r)=
\ep_0^{-\theta_0}
r^{-\al}>
(4M)^{\frac{1}{p-1}}r^{\frac{p-4}{p-1}}
\ ,$$
in view of
\eqref{eq-data-a=0-1} and $\al<(4-p)/(p-1)$.
To show $X_a$ is closed, it suffices to improve the lower bound. For that purpose, as 
$\pt W(0, r^*)>0$ and $G>0$,
we use D'Alembert's formula for \eqref{eq-Sch-1d-2} to obtain
\begin{align*}
W(t, r^*)\geq&\frac{W(0, r^*+t)+W(0, r^*-t)}{2}+\frac 12 \int^t_0\int^{r^*+(t-\tau)}_{r^*-(t-\tau)}
G(y, W)
 dy d\tau\\
\geq &
\frac 12 \ep_0^{-\theta_0}
(r^*-t)^{-\al}\\
\geq &
\frac 12\ep_0^{-\theta_0}
(r^*-t)^{\frac{p-4}{p-1}}
\gg (4M)^{\frac{1}{p-1}}r^{\frac{p-4}{p-1}}\ ,
\end{align*} where
we have used the fact that, 
 $r^*-t\le r^*\le 2r$
and $\ep_0\ll 1$. 
This proves $X_a=[0, L)\cap [0, a]$ for any $a<\infty$, and thus
we have \eqref{eq-Sch-lower-1} in $\Om$.

\subsection{Solution in the outgoing region: improved lower bound near the light cone}\label{sec-Sch-bu}
As we see from the proof of
\eqref{eq-Sch-lower-1}, we have actually proved a better lower bound for the solution:
\beeq\label{eq-Sch-lower-2}
W(t,r^*)>
\frac 12 \ep_0^{-\theta_0}
(r^*-t)^{-\al}
\ge \frac 12\ep_0^{-\theta_0} (r^*+t)^{-\al}
>0\ ,
\eneq
as well as the following partial differential inequality:
\beeq\label{eq-Sch-1d-3}
(\pa_t^2 -
\pa_{r^*}^2)W\ge \frac 13 (r^*)^{1-p} W^p,
\forall (t,r^*)\in\Om
\ .
\eneq
Equipped with the initial lower bound, \eqref{eq-Sch-lower-2}, we are ready to boost the bound to better bounds. 

Actually, under  the assumption
$$W(t, r^*)\ge C t^b(r^*+t)^{-a}, $$
with $a, b\ge 0$, we obtain,
\begin{align*}
W(t, r^*)\geq&\frac 12 \int^t_0\int^{r^*+(t-\tau)}_{r^*-(t-\tau)}
\frac 13 y^{1-p} W^p(\tau,y)
 dy d\tau\\
\geq &
\frac {C^p}6 
\int^t_0\int^{r^*+(t-\tau)}_{r^*-(t-\tau)}
(y+\tau)^{1-p(1+a)}\tau^{bp}
 dy d\tau\\
\geq &
\frac {C^p}3
(r^*+t)^{1-p(1+a)}
\int^t_0
\tau^{bp}(t-\tau)
 d\tau\\
= &
\frac {C^p}{3(bp+1)(bp+2)} 
(r^*+t)^{1-p(1+a)}
t^{bp+2}\ .
\end{align*} 
By induction, we see that
we could iterate
the initial lower bound, \eqref{eq-Sch-lower-2}, with 
$C=\ep_0^{-\theta_0}/2$, $a=\al$ and $b=0$, to the following bound, for any $m\in \N$,
\beeq\label{eq-Sch-lower-3}
W(t, r^*)\ge C(m)\ep_0^{-{\theta_0 p^m}}
(r^*+t)^{1-(\al+1)p^m} t^{2\frac{p^m-1}{p-1}} \ ,
 \forall (t,r^*)\in\Om\ .
\eneq
Here, we observe that $C(m)$ is
also dependent on $p$,  but is
 independent of $\al$ and $\ep_0$.

In particular, we observe that
the lower bound grows fast in the restricted region near the light cone,
$$\Om_1=\{
(t,r^*)\in\Om, r^*\le 2t\}
 \ ,$$
where we have $t\ge (t+r^*)/3$ and 
\beeq\label{eq-Sch-lower-4}
W(t, r^*)\ge C(m)\ep_0^{-{\theta_0 p^m}}3^{-2\frac{p^m-1}{p-1}}
(r^*+t)^{1-\frac{2}{p-1}+(\frac{2}{p-1}-1-\al)p^m} \ ,
 \forall (t,r^*)\in\Om_1\ .
\eneq

\subsection{Blow up near the light cone}\label{sec-2.4}
In view of the  lower bound, \eqref{eq-Sch-lower-4}, it is natural to expect that the solution will blow up near the light cone. 

To illustrate the blow up phenomenon for the solution, we try to extract ordinary differential inequality from the partial differential inequality \eqref{eq-Sch-1d-3}.

As we see the solution grows fast near the light cone, with parameter $
\nu$, it is naturally to introduce the functional, 
\beeq\label{eq-Sch-test}
G(t)=\int_0^1 \la^{\nu} W(t-\la, t+2\ep^{-N}+\la)d\la\ , t>1.\eneq
Then
\beeq\label{eq-Sch-test1}
G'(t)=\int_0^1 \la^{\nu} [(\pt+\pa_{r^*}) W](t-\la, t+2\ep^{-N}+\la)d\la\ .\eneq

Let
$Z=(\pt+\pa_{r^*}) W$, we have
\beeq\label{eq-Sch-1d-4}(\pt-\pa_{r^*}) Z=(\pt^2-\pa_{r^*}^2) W\ge \frac 13 (r^*)^{1-p} W^p \ ,\ \forall (t,r^*)\in\Om_1\ .
\eneq
To avoid the appearance of $r^*$, we use the improved bound
\eqref{eq-Sch-lower-4}, 
which tells us that there exists $m>0$ so that
$$1-\frac{2}{p-1}+(\frac{2}{p-1}-1-\al)p^m\ge 2\ .$$
For such $m$, we have
\beeq\label{eq-Sch-lower-5}
W(t, r^*)\ge C(m)\ep_0^{-{\theta_0 p^m}}3^{-2\frac{p^m-1}{p-1}}
(r^*)^{2} \ , \forall (t,r^*)\in\Om_1\ ,\eneq
and so
\beeq\label{eq-Sch-1d-5}
\frac 13 (r^*)^{1-p} W^p
\ge
\frac 13(C(m)\ep_0^{-{\theta_0 p^m}}3^{-2\frac{p^m-1}{p-1}})^{\frac{p-1}2}W^{\frac{p+1}2}\ge
\ep_0^{-\theta_0 \frac{p-1}2}W^{\frac{p+1}2}\ ,
\eneq
where, in the last inequality, we have assumed
$\ep_0\ll 1$ so that
$$\frac 13(C(m)\ep_0^{-{\theta_0 p^m}}3^{-2\frac{p^m-1}{p-1}})^{\frac{p-1}2} \ge
\ep_0^{-\theta_0 \frac{p-1}2}\ .$$
As we have observed, 
 $C(m)$ is 
 independent of $\al$, and all of the requirement on $\ep_0$ could be chosen to be independent of $\al$.
 
 With the help of
\eqref{eq-Sch-1d-4}
and 
\eqref{eq-Sch-1d-5}, we see that
\beeq\label{eq-Sch-1d-6}(\pt-\pa_{r^*}) Z\ge 
\ep_0^{-\theta_0 \frac{p-1}2}W^{\frac{p+1}2}
\ ,\ \forall (t,r^*)\in\Om_1\ ,
\eneq
and so, 
\beeq\label{eq-Sch-1d-7}
Z(t, 3A-t)-Z(A, 2A)\ge\int_A^t
\ep_0^{-\theta_0 \frac{p-1}2}W^{\frac{p+1}2}(\tau, 3A-\tau)d\tau
\ ,\ 
\eneq
for any $A\ge 2\ep^{-N}$ and $A\le t\le 3A/2-\ep^{-N}$.
In other words, for any $(t,r^*)\in \Omega_1$,
\beeq\label{eq-Sch-1d-8}
Z(t, r^*)-Z(\frac{t+r^*}3, 2\frac{t+r^*}3)\ge\int_{(t+r^*)/3}^t
\ep_0^{-\theta_0 \frac{p-1}2}W^{\frac{p+1}2}(\tau, t+r^*-\tau)d\tau
\ . 
\eneq

Concerning $Z=(\pt+\pa_{r^*}) W$, we recall that our condition on the initial data
\eqref{eq-data-a=0-1}
 has ensured that
$$
Z(0, r^*)>0, \forall r^*\ge 2\ep^{-N}\ ,
$$
provided that $\theta_1>\theta_0$ and $\ep_0\ll 1$.
By 
\eqref{eq-Sch-1d-4}, we see that
$$Z(t, r^*)\ge Z(0, r^*+t)>0, \forall (t, r^*)\in \Omega\ .$$
In particular, together with
\eqref{eq-Sch-1d-8}, we find that
\beeq\label{eq-Sch-1d-9}
Z(t, r^*)\ge \ep_0^{-\theta_0 \frac{p-1}2}\int_{0}^{t-(t+r^*)/3}
W^{\frac{p+1}2}(t-\tau, \tau+r^*)d\tau
\ ,\  \forall (t, r^*)\in \Omega_1\ .
\eneq

To ensure all of the $(t,r^*)$, appeared
in the definition of $G(t)$,  lies in $\Om_1$, we restrict ourselves to the case 
$t\ge 6+2\ep^{-N}$.
In particular,
for any $t\ge 6+2\ep^{-N}$ and $\la\in [0,1]$,
we have
$(t-\la)-((t-\la)+t+2\ep^{-N}+\la)/3\ge 1$ and so
\beeq\label{eq-Sch-1d-10}
Z(t-\la, t+2\ep^{-N}+\la)\ge \ep_0^{-\theta_0 (p-1)/2}\int_{0}^{1}
W^{(p+1)/2}(t-\la-\tau, t+2\ep^{-N}+\la+\tau)d\tau
\ .  
\eneq

Recalling
\eqref{eq-Sch-test1}, we get
\begin{eqnarray*}
G'(t) & = &\int_0^1 \la^{\nu} Z(t-\la, t+2\ep^{-N}+\la)d\la \\
 & \ge & \ep_0^{-\theta_0 (p-1)/2} \int_0^1 \la^{\nu} \int_{0}^{1}
W^{(p+1)/2}(t-\la-\tau, t+2\ep^{-N}+\la+\tau)d\tau d\la
\\
 & \ge & \ep_0^{-\theta_0 (p-1)/2} \int_{0}^{1} \la^{\nu} \int_\la^{\la+1}
W^{(p+1)/2}(t-y, t+2\ep^{-N}+y)dy d\la
\\
 & \ge & \ep_0^{-\theta_0 (p-1)/2} \int_0^{1}  \int_{0}^{y}
\la^{\nu}W^{(p+1)/2}(t-y, t+2\ep^{-N}+y) d\la dy
\\
 & \ge & \ep_0^{-\theta_0 (p-1)/2} \int_0^{1}  \frac{y^{\nu+1}}{\nu+1}
W^{(p+1)/2}(t-y, t+2\ep^{-N}+y)  dy\ .
\end{eqnarray*}

To connect the right hand side with $G$, we use H\"older's inequality to find
\begin{eqnarray*}G(t)
&= &\int_0^1  \la^{\frac{2}{p+1}} W(t-\la, t+2\ep^{-N}+\la)
\la^{-\frac{2}{p+1}}
\la^{\nu}
d\la\\
& \le&
\|\la^{\frac{2}{p+1}} W(t-\la, t+2\ep^{-N}+\la)\|_{L^{(p+1)/2}( \la^\nu d\la)}
\|\la^{-\frac{2}{p+1}}
\|_{L^{(p+1)/(p-1)}
( \la^\nu d\la)}\\
&\le&\|\la^{\frac{2}{p+1}} W(t-\la, t+2\ep^{-N}+\la)\|_{L^{(p+1)/2}( \la^\nu d\la)}
\ ,
\end{eqnarray*}
if we take $\nu=2/(p-1)$.

Combining these two estimates, we obtain
\beeq\label{eq-Sch-1d-11}
G'(t)\ge \frac{p-1}{p+1}\ep_0^{-\theta_0\frac{p-1}2} G^{\frac{p+1}2}\ ,\ 
\forall t\in [ 6+2\ep^{-N}, L)\ ,\eneq
which gives us
$$
\pa_t G^{-\frac{p-1}2}= -\frac{p-1}2 G^{-\frac{p+1}2}G'(t)\le 
-\frac{(p-1)^2}{2(p+1)}\ep_0^{-\theta_0 \frac{p-1}2}
\ ,\ 
\forall t\in [ 6+2\ep^{-N}, L).
$$
and
\beeq\label{eq-Sch-1d-12}
L\le 6+2\ep^{-N}+
\frac{2(p+1)}{(p-1)^2}\ep_0^{\theta_0\frac{p-1}2}
G^{-\frac{p-1}2}(6+2\ep^{-N})\ .
\eneq
 
 Recalling
\eqref{eq-Sch-test}
and \eqref{eq-Sch-lower-5}, it is clear that 
$$
 G(6+2\ep^{-N})\ge
  C(m)\ep_0^{-{\theta_0 p^m}}3^{-2\frac{p^m-1}{p-1}}
  \int_0^1 \la^{\nu} (6+4\ep^{-N})^2d\la  
  \ge \ep_0^{-\theta_0}\ep^{-2N}\ ,
 $$ 
 for sufficiently small $\ep_0$.
 
 In summary, we see that
$$ L\le 6+2\ep^{-N}+
\frac{2(p+1)}{(p-1)^2}\ep_0^{\theta_0(p-1)}
\ep^{N (p-1)}\le 3\ep^{-N}\ .$$
As this is achieved for any
$\ep, \mu\in (0,\ep_0)$, for some $\ep_0>0$. By taking limit $\mu\to 0$, we obtain
$$ T_*\le \lim_{\mu\to 0}L(\mu)\le  3\ep^{-N_0}\ ,$$
which completes the proof.
\begin{rem}
Due to this upper bound of the lifespan, 
a standard domain of dependence argument ensures that
we could choose the initial data of the form
\eqref{eq-data-a=0-1} in the region $r\in [\ep^{-N}, 10\ep^{-N}]$, and then extend as compactly supported smooth functions.
 \end{rem}
\subsection{The phenomenon of finite time blow-up for $1<p<p_c$, alternative approach}\label{sec-2.5}
It turns out that, to show blow up, we could avoid the ordinary differential inequality. Instead, we try to obtain a more precise lower bound for the constant
$C(m)
$, appeared in
\eqref{eq-Sch-lower-3}.

We know that
we could choose $C(0)=1/2$ and
$$C(m+1)=\frac{C(m)^p}{3(b_m p+1)(b_m p+2)}\ ,$$
where $b_m=2\frac{p^m-1}{p-1}$.

Let $d_0=2$, $d_1=2^{p+3}$, 
$d_m=d_1^{p^{m-1}}\Pi_{j=1}^{m-1}(4p b_j)^{2 p^{m-1-j}}$
for $m\ge 2$,
we observe that\footnote{
The case $m=0,1$ is trivial.
Suppose it is true for some  $m\ge 1$, then
$d_m^p=d_1^{p^{m}}\Pi_{j=1}^{m-1}(4p b_j)^{2 p^{m-j}}$.
Noticing that $b_m\ge 2$, we have
$b_m p+1\le b_m p+2\le b_m (p+1)\le 2b_m p$,
$$C(m+1)\ge \frac{1}{d_m^{p}
3(b_m p+1)(b_m p+2)
}
\ge
 \frac{1}{d_m^{p}
(4b_m p)^2}
=
 \frac{1}{d_1^{p^{m}}(4b_m p)^2\Pi_{j=1}^{m-1}(4p b_j)^{2 p^{m-j}}
}=\frac{1}{d_{m+1}}\ .
$$}
$$C(m)\ge \frac{1}{d_m}\ .$$
Recall  that
$$b_m=2\frac{p^m-1}{p-1}\le \frac{2}{p-1}p^m\ ,
$$
we have
 \begin{eqnarray*}
\ln d_m & = & p^{m-1}((p+3)\ln 2+\Sigma_{j=1}^{m-1} 2{p^{-j}}\ln (4pb_j))\\
 & \le & p^{m-1}((p+3)\ln 2+\Sigma_{j=1}^{m-1} 2{p^{-j}}(\ln 
\frac{8p}{p-1}+j\ln p)
)\ ,
\end{eqnarray*}
from which
we see that there exists $C=C(p)$ so that
$$\ln d_m\le C(p) p^m\ .$$ In conclusion, we obtain a lower bound of $C(m)$: $$C(m)\ge \exp(-C(p)p^m)\ .$$

If $\Om_1\neq\emptyset$, by \eqref{eq-Sch-lower-4}, we have
for $(t, r^*)=(2\ep^{-N},4\ep^{-N})$ and any $m\in\N$,
\beeq\label{eq-Sch-lower-6}
W(2\ep^{-N},4\ep^{-N})\ge 
(6\ep^{-N})^{1-\frac{2}{p-1}}
\exp({p^m((\frac{2}{p-1}-1-\al)\ln (\ep^{-N})+\ln(3^{-\frac{2}{p-1}}\ep_0^{- \theta_0})-C(p))}) \ .
\eneq

For $\ep_0\ll 1$, independent of $\al\le \frac{2}{p-1}-1$,
 we have
$$(\frac{2}{p-1}-1-\al)\ln (\ep^{-N})+\ln(3^{-\frac{2}{p-1}}\ep_0^{- \theta_0})-C(p)
\ge
\ln(3^{-\frac{2}{p-1}}\ep_0^{- \theta_0})-C(p)
>0\ .$$
Letting $m\to \infty$, we see that
  the solution blows up at
$(t, r^*)=(2\ep^{-N},4\ep^{-N})$, which shows that
$\Omega_1=\emptyset$. In particular, we have
$$T_*\le T_*(u_0, u_1)\le 2\ep^{-N}, \forall \ep, \mu\in (0,\ep_0)\ . $$
Let $\mu\to 0$, we obtain
$ T_*\le  2\ep^{-N_0}$ as desired.

\section{Kerr space-time}\label{sec-Ker}
\subsection{Kerr space-time}\label{sec-Ker1}

The Kerr metric in the Boyer-Lindquist coordinates takes the  form 
\eqref{metric-Kerr}.
For future reference, we record the nontrivial components of the inverse of $g_{\al\be}(r,\theta)$: 
$$g^{tt}=-\frac{(r^2+a^2)^2-a^2\La\sin^2\theta}{\rho^2\La}\ , \ g^{t\phi}=-a\frac{2Mr}{\La\rho^2}\ ,$$
$$g^{rr}=\frac{\La}{\rho^2}\ , \ g^{\theta\theta}=\frac{1}{\rho^2}\ , \ g^{\phi\phi}=\frac{\La-a^2\sin^2\theta}{\rho^2\La\sin^{2}\theta}\ .$$
We observe that
we have \eqref{eq-stru-Kerr} for the Kerr metric, as it is stationary $\pt g_{\al\be}=0$ and $g_{\al\be}$ is independent of $\phi$.

\subsection{Initial data }
Let $0<\ep\le \ep_0$, $0<\theta_0<\theta_1<1$, and $L=\ep^{-\frac{p(p-1)}{1+2p-p^2}}$. We choose the initial data taking the form  
\beeq\label{eq-Kerr-data}
u_0(r)=\ep_0^{-\theta_0}
r^{-\frac{2}{p-1}},  \ 
u_1(r)=\ep_0^{-\theta_1}
r^{-\frac{2}{p-1}-1}, \ 6L\leq r\leq 8L\ .
\eneq
As before, it is easy to check that, provided that $\ep_0\ll 1$, there exists $C^\infty_0$ extension
so that
$$\supp\ u_0, u_1\subset \{r\in [5L, 10L]\},\ 
\|(u_0, u_1)\|_{X}\ \le \ \ep_0^{-1}\ep\  .$$

Thanks to the well-posedness
 and finite speed of propagation (with speed less than two), we infer that, for any $t\in [0, T_*(u_0, u_1))$, the solution $u(t)$ is compactly supported and
\beeq\label{eq-Kerr-spt}\supp\ u(t)\subset \{(r,\omega): r\in [5L-2t, 10L+2t]\}\ .\eneq

We claim $T_*(u_0, u_1)\le L$, for $\ep_0\ll 1$.
Suppose by contradiction, that $T_*(u_0, u_1)> L$, then for any
$t\in [0, L]$, we have
\beeq\label{eq-Kerr-spt2}\supp\ u(t)\subset \{(r,\omega): r\in [5L-2t, 10L+2t]\} 
\subset \{(r,\omega): r\in
 [3L, 12L]\}\ .\eneq
In addition,
$$u\in C([0, L]; H^1)\cap
C^1([0, L]; L^2)
\cap
C^2([0, L]; H^{-1})
\ .
$$

\subsection{Blow-up}\label{sec:Ker-bu}

Recall that
with $F_p(u)=|u|^p$,
the Cauchy problem \eqref{nlw} is in the following form
$$g^{-1/2}\pa_\al g^{1/2}g^{\al\be}\pa_\be u=-|u|^p\ .$$
Similar to subsection \ref{sec-2.4}, we try to extract ordinary differential inequality from the PDE.

For that purpose, we try to exploit the support information
\eqref{eq-Kerr-spt2} and
use integration by parts. To be more specific, by integration in space, with respect to the natural volume form $g^{1/2}dx$, we obtain
\begin{eqnarray*}
\int_{r\in [3L, 12L]} |u|^p g^{1/2}dx & = & 
-\int_{r\in [3L, 12L]} \pa_\al (g^{1/2}g^{\al\be}\pa_\be u) dx \\
  & = &-\pa_t \int_{r\in [3L, 12L]} g^{1/2}g^{t\be}\pa_\be u dx\ .
  \end{eqnarray*}
As we have observed in Subsection \ref{sec-Ker1}, we have
\eqref{eq-stru-Kerr}, which implies
$$\pa_\be(g^{1/2}g^{t\be} u)=
g^{1/2}g^{t\be} \pa_\be u\ .$$
Thus, we conclude the following ordinary differential equation:
$$\int_{r\in [3L, 12L]} |u|^p g^{1/2}dx=
\pa_t^2 \int_{r\in [3L, 12L]} (-g^{tt}) u g^{1/2} dx\ .$$

 Let $W(t)=
 \int_{r\in [3L, 12L]} (-g^{tt}) u(t) g^{1/2} dx\in C^2$. By H\"older's inequality, we get 
 $$W(t)\le \Big(\int_{r\in [3L, 12L]} |u|^p g^{1/2} dx\Big)^{\frac{1}{p}}\Big(\int_{r\in [3L, 12L]} g^{1/2} dx\Big)^{\frac{1}{p'}} \le C \Big(W''(t)\Big)^{\frac{1}{p}}L^{3/p'}\ ,$$
 for some universal constant $C>0$.
 That is,  we get an ordinary differential inequality:
\beeq\label{eq-Ker-odi}W'' \ge C^{-p}L^{-3(p-1)}W^p, \forall t\in [0, L]\ .\eneq

Recall the initial dat \eqref{eq-Kerr-data}, for
$\ep_0\ll 1$, we have
$-g^{tt}g^{1/2}\ge 1/2$ and
there exists some universal constant $c>0$ so that
\beeq
\label{eq-Kerr-data2}
W(0)=\int_{r\in [3L, 12L]} (-g^{tt}) u_0 g^{1/2} dx
\in ( c \ep_0^{-\theta_0} L^{3-\frac{2}{p-1}}, 
 c^{-1} \ep_0^{-\theta_0} L^{3-\frac{2}{p-1}})\ ,
\eneq 
and
\beeq
\label{eq-Kerr-data3}
W'(0)
\ge c \ep_0^{-\theta_1} L^{2-\frac{2}{p-1}}\ .
\eneq

Based on \eqref{eq-Ker-odi}, and
\eqref{eq-Kerr-data2}-\eqref{eq-Kerr-data3}, it is clear that
$W, W'>0$ for any $t\in [0, L]$.
Multiplying $F'$, we see that
\beeq\label{eq-Ker-odi2}
\frac{(W')^2}2-C^{-p}L^{-3(p-1)}\frac{W^{p+1}}{p+1}
\ge 
\frac{(W')^2(0)}2-C^{-p}L^{-3(p-1)}\frac{W^{p+1}(0)}{p+1}>0
, 
\eneq
provided that
$2\theta_1>(p+1)\theta_0$ and $\ep_0\ll 1$.

By \eqref{eq-Ker-odi2}, we have
\beeq\label{eq-Ker-odi3}
W'>\sqrt{\frac{2}{(p+1)C^{p}}}L^{-\frac{3(p-1)}2}W^{\frac{p+1}2},\ 
 \forall t\in [0, L]\  ,
\eneq
which tells us that
$$\sqrt{\frac{2}{(p+1)C^{p}}}L^{-\frac{3(p-1)}2}t\le
\frac{2}{p-1} W^{-\frac{p-1}2}(0)
\le
\frac{2}{p-1} c^{-\frac{p-1}2}
\ep_0^{\frac{p-1}2\theta_0} L^{1-\frac{3(p-1)}2}\ ,
$$
for any $t\in [0, L]$.
Taking $t=L$, we get
$$\sqrt{\frac{2}{(p+1)C^{p}}}\le
\frac{2}{p-1} c^{-\frac{p-1}2}
\ep_0^{\frac{p-1}2\theta_0} \ ,
$$
which is the desired contradiction, for $\ep_0\ll 1$.
This completes the proof of
Theorem \ref{thm-Kerr}.

\section{Asymptotically flat manifolds}\label{sec-AF}
In this section, we give the proof of
Theorem \ref{thm-AF} with $F_p(u)=|u|^p$, by adapting the proof of
Theorem \ref{thm-Kerr} to general $(1+n)$ dimensional asymptotically flat manifolds, with
structural condition 
\eqref{eq-stru-Kerr}.

\subsection{Initial data and preparation}
Let $n\ge 2$, $\theta_0, \theta_1\in (0,1)$ to be specified later
and $$L=\ep^{-\frac{2p(p-1)}{2+(n+1)p-(n-1) p^2}}\ .$$
For $r\in [6L, 8L]$, we set the data as for Kerr:
\eqref{eq-Kerr-data}.
Observing that, for any $m\ge 1$, $k\geq 0$,
$$
\sum_{|\ga|\le k} \| (\nabla, \Omega)^{\ga}
 (\nabla u_0^\ep, u^\ep_1)\|_{\dot H^{m-1}(r\in [6L,8L])}\sim L^{-m-2/(p-1)+n/2}$$
and $L^{n/2-2/(p-1)-s_d}=\ep$.
It is clear that, for any $k\ge 0$, there exists $\ep_0\ll 1$ such that
the initial data
have 
extension in $C^\infty_0(r\in (5L, 10 L))$ so that
\beeq
\label{eq-initial2}
\sum_{|\ga|\le k} \| (\nabla, \Omega)^{\ga}(\nabla u_0, u_1)\|_{L^2\cap\dot{H}^{s_d-1}}\le \ep_0^{-1} \ep\ , 
\ s_d=\frac{1}{2}-\frac{1}{p}\ .
\eneq
We will fix such choice of data.

In view of the asymptotically flat assumption
\eqref{eq-H1}, we know that the maximum speed of propagation for the $\Box_\gm$ is less than $2$, in the region $r\ge L$, provided that $\ep_0\ll 1$.
With the help of this fact, we know from finite speed of propagation for the energy solution,
$u=u(t,x)$, that
$$\supp\ u(t,\cdot)\subset \{x: r\in (5L-2t, 10L+2t)\}\subset \{x: r> 3 L\}, \forall t\in [0,L]\cap [0, T_*)\ .$$
Moreover, we have
$L>R$, such that
$$M^t_R\supset \{r\ge L+2t\}
\supset \{r> 3L\}, \forall t\in [0,L]\ .
$$

\subsection{Ordinary differential inequalities}
To illustrate the result $T_*\le L$, we proceed with contradiction, by assuming $T_*>L$.
Then we have $$u\in C([0,L]; H^1(r>3L))\cap C^1([0,L]; L^2(r>3L))\cap 
C^2([0,L]; H^{-1}(r>3L))\ .$$

With the help of the equation,
we see that
\begin{eqnarray*}
-\int_{r>3L} F_p(u) g^{1/2}dx & = & 
\int_{r>3L} (\Box_{\gm} u) g^{1/2}dx 
 \\
  & = & \pa_t\int_{r>3L} 
[\pa_\be(g^{1/2}g^{0\be} u)-
\pa_\be(g^{1/2}g^{0\be}) u]
 dx \\
  & = & \pa_t^2\int_{r>3L}  g^{1/2}g^{00} u dx \ ,
\end{eqnarray*}
where we have used  the structural condition  \eqref{eq-stru-Kerr} in the last equality.
In view of the differential identity,
 we introduce the auxillary functional
$$W(t):=-\int g^{1/2} g^{00}u(t,x) dx\in C^2([0,L])\ ,
$$
and rewrite the identity as follows
$$W''(t)=\int F_p(u)g^{1/2}dx
= \|u\|^p_{L^p (g^{1/2}dx)}\ .$$
In addition, for any given $\nu>0$, 
H\"older's inequality gives us
$$W(t)\le
\|u(t)\|_{L^p(g^{1/2}dx)}\ep_0^{-\nu/p} L^{n/p'}
\ ,$$
for any $\ep_0\ll 1$.

In summary, we see that
\beeq\label{eq-AF-ODI}
 W''(t)
\ge \|u\|^p_{L^p (g^{1/2}dx)}
\ge \ep_0^\nu W^pL^{-n(p-1)}\ .
\eneq

Concerning the data, in view of \eqref{eq-Kerr-data} and \eqref{eq-H1}, we have
\beeq \label{eq-AF-data} 
W(0)\sim
 \ep_0^{-\theta_0}L^{n-\frac{2}{p-1}},
\eneq
\beeq\label{eq-AF-data-2} 
W'(0)=
-\int (g^{1/2} g^{00}u_t(0,x) +(g^{1/2} g^{00})_t(0)u(0,x)) dx
\ge
 \ep_0^{-\theta_1/2}L^{n-\frac{2}{p-1}-1}\ ,
\eneq
provided that $\theta_1\ge \theta_0>0$ and $\ep_0\ll 1$.

\subsection{Blow-up}
Based on \eqref{eq-AF-ODI}, \eqref{eq-AF-data}-\eqref{eq-AF-data-2}, with appropriate choices of the parameters, we will illustrate the contradiction,
as in Subsection \ref{sec:Ker-bu}.

At first,  $W, W'>0$ for any $t\in [0, L]$. Then \eqref{eq-AF-ODI} tells us that
\begin{eqnarray*}
\frac{(W')^2}2- 
 \ep_0^{\nu} \frac{W^{p+1}}{(p+1)L^{n(p-1)}}
&\ge&
\frac{(W')^2(0)}2- 
 \ep_0^{\nu} \frac{W^{p+1}(0)}{(p+1)L^{n(p-1)}}\\
&\ge&
\left(\frac{ \ep_0^{-\theta_1}}2- C
 \ep_0^{\nu} \frac{ \ep_0^{-\theta_0 (p+1)}}{p+1}\right)L^{2(n-1-2/(p-1))}\ ,
\end{eqnarray*}
for some constant $C$, depending only on $n$.
The right hand side is positive, if $ \ep_0\ll 1$ and \beeq\label{eq-AF-choice1}
\theta_1>\theta_0 (p+1)-\nu\ .\eneq

As a consequence, we see that
\beeq\label{eq-AF-ODI2}
W'\ge
 \ep_0^{\nu}  W^{(p+1)/2} L^{-n(p-1)/2},
W(0)\ge 
 \ep_0^{-2\theta_0/3}L^{n-2/(p-1)}\ ,
\eneq
which gives us
$$
 \ep_0^{\nu}  t L^{-n(p-1)/2}\le
\frac{W(0)^{-(p-1)/2}-W(t)^{-(p-1)/2}}{(p-1)/2}\le 
\frac{2}{p-1} \ep_0^{(p-1)\theta_0/3}L^{1-n(p-1)/2},
$$
for any $t\in [0, L]$.
In conclusion, taking $t=L$, we obtain 
$$\frac{p-1}{2}\le 
 \ep_0^{(p-1)\theta_0/3-\nu}, \forall \ep\in (0, \ep_0),  \ep_0\ll 1\ ,$$ which is impossible if we have
\beeq\label{eq-AF-choice2}
(p-1)\theta_0/3-\nu>0\ .\eneq

A sample choice of the parameters for \eqref{eq-AF-choice1} and \eqref{eq-AF-choice2} could be
$$
\theta_0=1/(p+2),\ \theta_1=(p+1)\theta_0,\ 
\nu=\frac{p-1}{4}\theta_0\ .
$$
This completes the proof.

 \section{
Spherically symmetric asymptotically flat manifolds 
 }\label{sec-AF-radial}
 Inspired by the proof for the Schwarzschild space-time, we intend to generalize it to spherically symmetric asymptotically flat manifolds.
 
\subsection{Initial data and transformation}

Let $n\ge 2$, $p\in (1,p_c)$ with
$$\frac{1}{p}<\al<\frac{2}{p-1}-\frac{n-1}{2}\ , \ N=\frac{p}{\al p-1}\ . $$
As $p\in (1,p_c)$, we always have
$$
\frac{2}{p-1}-\frac{n-1}{2}>\frac 1p\ .$$

As above,  for the part near spatial infinity, $r\geq \ep^{-N}$, we choose the outgoing initial data 
\beeq\label{eq-data-radi}
u_0=
 \ep_0^{-\theta_0}
r^{-\al-\frac{n-1}{2}}, \ u_1= \ep_0^{-\theta_1} r^{-\al-1-\frac{n-1}{2}}\ ,
\eneq
where $0<2\theta_0<\theta_1<1$, and $\ep_0>0$ is a small parameter to be determined later.

In the polar coordinates, we assume $\gm$ is of the form
\eqref{metric-sph}, that is,
$$\gm=g_{tt}(t,r)dt^2+g_{rr}(t,r)dr^2+2g_{tr}(t,r)dtdr+ r^2d\omega^2, \omega\in \Sp^{n-1}\ .$$
Let $G=\sqrt{g_{tr}^2-g_{tt}g_{rr}}$, we have
$$g^{1/2}=|\det \gm|^{1/2}=Gr^{n-1},
g^{tt}=g_{rr}/G^2,
g^{tr}=-g_{tr}/G^2,
g^{rr}=g_{tt}/G^2
\ .$$

Then $\Box_\gm u$ with radial $u(t,r)$ reads 
\beeq
\label{ksl}\Box_\gm u=
g^{-1/2}\pa_\al (g^{1/2}g^{\al\be}\pa_\be u)=
g^{tt}\pa^{2}_t u+2g^{tr} u_{tr}
+g^{rr}\pa^{2}_r u
+C^tu_t+C^r u_r\ ,\eneq
where 
$$
C^\al(t, r)=\pa_\be g^{\al\be}+ g^{\al\be}\frac{\pa_\be G}{G}
+ g^{\al r}\frac{n-1}{r}\ .
$$
Due to the assumption \eqref{eq-H2}, we have
$$C^t=\CO(r^{-1-\de_0}), C^r=g^{r r}\frac{n-1}{r}+
\tilde C^r=g^{r r}\frac{n-1}{r}+
\CO(r^{-1-\de_0})\ .$$
Recalling that
$$\pa_r^2+\frac{n-1}{r}\pa_r=r^{-(n-1)/2}\pa_r^2 r^{(n-1)/2}-\frac{(n-1)(n-3)}{4r^2}, $$
we find that, with $W=r^{(n-1)/2}u$,
\begin{eqnarray*}
r^{(n-1)/2}\Box_g u&=&
r^{(n-1)/2}\Box_g (r^{-(n-1)/2} W)\\
&=&g^{tt}\pa^{2}_t W+2g^{tr} (W_{tr}
-\frac{n-1}{2r}W_t)
+g^{rr}(W_{rr}-\frac{(n-1)(n-3)}{4r^2}W)\\&&
+C^t W_t+\tilde C^r(t, r)
(W_r-\frac{n-1}{2r}W)\\
&=&
g^{tt}W_{tt}+g^{rr}W_{rr}+2g^{tr}W_{tr}+K_1W_r+K_2W_t+K_3W
\ ,\end{eqnarray*}
with
$$K_1, K_2=\CO( r^{-1-\de_0})\ , K_3=\CO( r^{-2})\ , \forall r\gg 1\ .$$

Thus, the equation
$\Box_{g}u+F_p(u)=0$ is equivalent to the following
\beeq\label{eq-afrad-1} g^{tt}W_{tt}+g^{rr}W_{rr}+2g^{tr}W_{tr}+K_1W_r+K_2W_t+K_3W=
-r^{\frac{n-1}2}F_p(r^{-\frac{n-1}2}W)\ .
\eneq

\subsection{Double null frame}
Introducing the double null frame, the equation for $W$ will be simplified and essentially solvable.

Let $\la_{+}(t,r)>\la_{-}(t,r)$ be the roots of
the characteristic equation $g^{tt}\la^2+2g^{tr}\la+g^{rr}=0$,
and $\eta, \xi$ be the corresponding ingoing/outgoing null coordinates, so that
$$\eta_t-\la_-(t,r)\eta_r=0=\xi_t-\la_+(t,r)\xi_r\ ,\ 
\eta(0,r)=-r, \xi(0,r)=r\ .
$$
Thanks to the assumption \eqref{eq-H2}, 
we have $\la_\pm=\pm 1+\CO(r^{-\de_0})$,
and it is then clear that 
$\eta_t, -\eta_r, \xi_t, \xi_r\sim 1$,
which tells us that
$\eta, \xi$ behaves like $t-r, t+r$.

With the help of the coordinates $(\eta,\xi)$,
we have that 
$$\pt=\eta_t \pa_\eta+\xi_t\pa_\xi
=\la_-\eta_r \pa_\eta+\la_+ \xi_r\pa_\xi
\ , \pa_r=\eta_r\pa_\eta+\xi_r\pa_\xi\ ,$$
$$
\pa_t-\la_+\pa_r=(\la_--\la_+)
\eta_r\pa_\eta, 
\pa_t-\la_-\pa_r=
 (\la_+-\la_-)\xi_r \pa_\xi\ .
$$
As 
$g^{tt}\la^2+2g^{tr}\la+g^{rr}=g^{tt}(\la-\la_+)(\la-\la_-)$, it is clear that
$$g^{tt}W_{tt}+g^{rr}W_{rr}+2g^{tr}W_{tr}
=g^{tt}(\pt-\la_+\pa_r)
(\pt-\la_-\pa_r)W+g^{tt}((\pt-\la_+\pa_r)\la_-)\pa_r W\ .
$$
Thus, 
in the outgoing region
\beeq\label{eq-afrad-region}\mathcal{K}:=\{(t,r): t\ge 0, \eta\le -\ep^{-N}\}
=\{(\xi,\eta): \xi+\eta\ge 0, \eta\le -\ep^{-N}\}
\ ,\eneq
the equation \eqref{eq-afrad-1} could be rewritten as
\beeq
\label{eq-afrad-2} \pa_\eta\pa_\xi W+C_1\pa_\eta W+C_2
\pa_\xi W
+C_3 W
=
\frac{ r^{(n-1)/2}}{
g^{tt}\eta_r\xi_r(\la_+-\la_-)^2}
F_p(r^{-\frac{n-1}2}W )\ ,
\eneq
with $C_1, C_2=\CO(\xi^{-1-\de_0})$ and
 $C_3=\CO(\xi^{-2})$.

To simplify the equation further, we introduce the integrating factor
\beeq
\label{eq-int-fac-K0}K(\xi, \eta)=e^{-\int^{-\ep^{-N}}_{\eta} C_2(\xi, \tau)d\tau}=1+\CO(\ep^{N\de_0})\ ,\eneq so that
$\pa_\eta K=C_2(\xi,\eta) K$ and $$\pa_\eta (KW_\xi)=K(W_{\xi\eta}+C_2 W_\xi)\ .$$
Multiplying $K$ to \eqref{eq-afrad-2} , we get
$$\pa_\eta (KW_\xi)+C_1 K\pa_\eta W
+C_3 K W
=
\frac{K  r^{(n-1)/2}}{
g^{tt}\eta_r\xi_r(\la_+-\la_-)^2}
F_p(r^{-\frac{n-1}2}W )\ ,$$
that is,
$$\pa_\eta (K(W_\xi +C_1 W))
+(C_3 K-\pa_\eta(C_1 K)) W
=
\frac{K  r^{(n-1)/2}}{
g^{tt}\eta_r\xi_r(\la_+-\la_-)^2}
F_p(r^{-\frac{n-1}2}W )\ .$$
Similarly,
with 
\beeq\label{eq-int-fac-K1}K_1(\xi, \eta)=e^{\int^{\xi}_{-\eta} C_1(\tau, \eta)d\tau}=1+\CO(\ep^{N\de_0})\ ,\eneq so that
$$(K_1W)_\xi=K_1(W_\xi+C_1 W)$$
and we have
$$\pa_\eta (\frac{K}{K_1}\pa_\xi (
K_1 W))
+(C_3 K-\pa_\eta(C_1 K)) W
=
\frac{K r^{(n-1)/2}}{
g^{tt}\eta_r\xi_r(\la_+-\la_-)^2}
F_p(r^{-\frac{n-1}2}W )\ .$$

Let $U=K_1 W$, we finally
obtain
the equation we use to do the integration, for which we write it in the form akin to 
\eqref{eq-Sch-1d-2},
\beeq
\label{z5'f}
\pa_\eta (\frac{K}{K_1}\pa_\xi U)
=
\frac{K}{
g^{tt}\eta_r\xi_r(\la_+-\la_-)^2}
 r^{(n-1)/2}
F_p(\frac{U}{K_1 r^{(n-1)/2}})
-\frac{C_3 K-\pa_\eta(C_1 K)}{K_1} U\ .
\eneq

Observing that
$\xi_r, -\eta_r=1+\CO(\ep^{N\de_0})$, we get
\beeq
\label{eq-afrad-3} 
\pa_\eta \left(\frac{K}{K_1}\pa_\xi U\right)
=
\frac{1+\CO(\ep^{N\de_0})}{
4}
 r^{(n-1)/2}
F_p\left(\frac{U}{K_1 r^{(n-1)/2}}\right)
+\CO(r^{-2}) U\ .
\eneq


\subsection{Solution in the outgoing region: initial lower bound}

Recalling our choice of the initial data \eqref{eq-data-radi}, by assuming
$\theta_1>2\theta_0>0$ and sufficiently small $\ep_0$,
we have
\beeq
\label{eq-data-radi1}
U(\xi,-\xi)\ge
 \ep_0^{-\theta_0/2}
\xi^{-\al}, \ U_\eta(\xi,-\xi)\ge U_\xi(\xi,-\xi)\ge 
 \ep_0^{-\theta_1/2} \xi^{-\al-1}\ .
\eneq

At first, we claim that, as long as
$\al<
2/(p-1)-(n-1)/2$,
\beeq\label{eq-induc1}
U(\xi,\eta)\ge  \ep_0^{-\theta_0/2}\xi^{-\al}, 
 \forall (\xi,\eta)\in\mathcal K\ ,
\eneq
where $\mathcal K$ is the region
\eqref{eq-afrad-region}.
By continuity, it suffices to prove
\eqref{eq-induc1}, in $K_{\xi_0}$, under the weaker
assumption
\beeq\label{eq-induc2}
U(\xi,\eta)>  \ep_0^{-\theta_0/3}\xi^{(n-1)/2-2/(p-1)}, \forall (\xi,\eta)\in\mathcal K_{\xi_0}\ ,
\eneq
where $$\mathcal K_{\xi_0}=\{(\xi,\eta): \ep^{-N}\le -\eta\le \xi\le \xi_0\}\ .$$

Actually,
by \eqref{eq-induc2}, we have $u>0$
in $\mathcal K_{\xi_0}$,
 and
so
 $F_p(u)=u^{p}$. Then
$$
r^{(n-1)/2}
F_p\left(\frac{U}{K_1 r^{(n-1)/2}}\right)\sim
r^{-(n-1)(p-1)/2}
U^p\gg r^{-2} U\ ,
$$
where we have used the fact that $r\le \xi\les r$ in
$\mathcal K$.
Thus,
by \eqref{eq-afrad-3},
 \beeq
\label{eq-afrad-4}
\pa_\eta \left(\frac{K}{K_1}\pa_\xi U\right)
\ge
\frac{U^p}{5r^{(n-1)(p-1)/2}} 
\ge
\frac{U^p}{5\xi^{(n-1)(p-1)/2}} 
\ ,\eneq
provided that $ \ep_0\ll 1$.
With the help of
\eqref{eq-afrad-4} and
\eqref{eq-data-radi1}, we could integrate from $(\xi, -\xi)$ to
$(\xi, \eta)$, to conclude
$$
\left.\left(\frac{K}{K_1}\pa_\xi U\right)(\xi, \la)\right|_{\la=-\xi}^\eta=\int_{-\xi}^\eta
\frac{U^p(\xi,\la)}{5\xi^{(n-1)(p-1)/2}} 
 d\la>0\ .
$$
Recalling that \eqref{eq-int-fac-K0}-\eqref{eq-int-fac-K1} ensures
$K/K_1\in (5/6, 6/5)$ and so
\beeq
\pa_\xi U(\xi, \eta)\ge 
\frac 12\pa_\xi U(\xi, -\xi)\ge
\frac 12\ep_0^{-\theta_1/2} \xi^{-\al-1}>0\ ,
 \forall (\xi,\eta)\in\mathcal K\ .
\eneq

By
integrating from $( -\eta, \eta)$ to
$(\xi, \eta)$, we obtain
$$
U(\xi, \eta)
\ge
U(-\eta, \eta)\ge \ep_0^{-\theta_0/2}(-\eta)^{-\al}
\ge \ep_0^{-\theta_0/2}\xi^{-\al}\ ,
$$
which completes the proof of 
\eqref{eq-induc1}.

\subsection{Solution in the outgoing region: improved lower bound near the light cone}
Equipped with the initial lower bound, \eqref{eq-induc1}, we are ready to boost the bound to better bound.

Under  the assumption
$$U\ge C \xi^{-a}(\xi+\eta)^b, $$
with $a, b\ge 0$, we obtain from
\eqref{eq-afrad-4} and
\eqref{eq-data-radi1} that
$$\frac{K}{K_1}\pa_\xi U(\xi,\eta)>\int_{-\xi}^\eta
\frac{U^p(\xi,\la)}{5\xi^{(n-1)(p-1)/2}} 
 d\la
\ge \frac{C^{p}(\xi+\eta)^{bp+1}}{5(bp+1)\xi^{ap+(n-1)(p-1)/2}} \ .
$$
As
$K/K_1\le 6/5$,
we could integrate further with respect to $\xi$:
$$U(\xi,\eta)>\int_{-\eta}^\xi
 \frac{C^{p}(\la+\eta)^{bp+1}}{6(bp+1)\la^{ap+(n-1)(p-1)/2}} 
 d\la\ge
 \frac{C^{p} (\xi+\eta)^{bp+2}}{6(bp+1)(bp+2)\xi^{ap+(n-1)(p-1)/2}} \ .
$$
By induction, we see that
we could iterate
the initial lower bound, \eqref{eq-induc1}, to the following improved bound, for any $m\in \N$,
\beeq\label{eq-induc3}
U(\xi,\eta)\ge C(m)\ep_0^{-\frac{\theta_0 p^m}2}\xi^{\frac{n-1}2-(\al+\frac{n-1}2)p^m} (\xi+\eta)^{2\frac{p^m-1}{p-1}} \ ,
 \forall (\xi,\eta)\in\mathcal K\ .
\eneq
Here, as in Subsection \ref{sec-Sch-bu}, $C(m)$ is
 independent of $\al$ and $\ep_0$.

In particular, in the restricted region,
$$\mathcal K_1=\{
(\xi,\eta)\in\mathcal K, \xi\ge -2\eta\}
 \ ,$$
we have $\xi+\eta\ge \xi/2$ and so
\beeq\label{eq-induc4}
U(\xi,\eta)\ge C(m)
(2^{-\frac{2}{p-1}}\ep_0^{-\frac{\theta_0 }2})^{p^m}
\xi^{\frac{n-1}2-\frac{2}{p-1}+(\frac{2}{p-1}-\frac{n-1}2-\al)p^m} \ ,
 \forall (\xi,\eta)\in\mathcal K_1\ ,
\eneq
for some $C(m)<\infty$.

\subsection{Blow up near the light cone}
With the help of the improved lower bounds near the light cone, 
\eqref{eq-induc4}, the similar proof as for the Schwarzschild space-time could be applied to show the blow up and the upper bound of the lifespan. And we shall omit the details.

\end{document}